\documentclass{amsart}
\usepackage{amsmath,latexsym,amssymb}

\newtheorem{thm}{Theorem}[section]
\newtheorem{lem}[thm]{Lemma}
\newtheorem{prop}[thm]{Proposition}
\newtheorem{rem}[thm]{Remark}
\newtheorem{cor}[thm]{Corollary}
\newtheorem{exam}[thm]{Example}

\begin{document}

\title{Spectrum of partial integral operators with degenerate kernel}

\author{Eshkabilov Yu.Kh.}
\address{Department of mechanics and mathematics. National University of Uzbekistan, Vuzgorodok, 100174, Tashkent, Uzbekistan.}
\email{yusup62@mail.ru}
\author{Arzikulov G.P.}
\address{Department of energetics, Tashkent State Technical University, Vuzgorodok, 100174, Tashkent, Uzbekistan.}
\email{arzikulov79@mail.ru}
\author{Haydarov F.H.}
\address{Department of mechanics and mathematics. National University of Uzbekistan, Vuzgorodok, 100174, Tashkent, Uzbekistan.}
\email{haydarov\_imc@mail.ru}
\newpage

\begin{abstract} In the paper we consider self-adjoint partial integral operators of Fredholm type $T$ with
a degenerate kernel on the space $L_2([a,b]\times[c,d]).$
Essential and discrete spectra of $T$ are described.
\end{abstract}

\keywords{spectrum, essential spectrum, discrete spectrum, partial
integral operator, partial integral equation.}

\subjclass[2010]{Primary: 45A05, 47A10, 47G10; Secondary: 45P05,
45B05, 45C05.}

\maketitle

\section{Introduction}\label{EshkArz_Section1}

Linear equations and operators involving partial integrals appear
in elasticity theory \cite{kal1,kov1,vek1}, continuum mechanics
\cite{kov1,kov2,kov3,kov4}, aerodynamics \cite{kal2} and in PDE
theory \cite{gour,mun}. Self-adjoint partial integral operators
arise in the theory of Schr\"{o}dinger operators
\cite{esh1,lak1,ras1,esh2}. Spectrum of a discrete Schr\"{o}dinger
operator $H$ are tightly connected (see \cite{esh2,esh3}) with
that of partial integral operators which participate in the
presentation of the operator $H.$

Let $\Omega_1$ and $\Omega_2$ be closed boundary subsets in
$\mathbb{R}^{\nu_1}$ and $\mathbb{R}^{\nu_2},$ respectively.
Partial integral operator (PIO) of Fredholm type in the space
$L_p(\Omega_1\times\Omega_2),$ $p\geq1$ is an operator of the form
\cite{app1}:
\begin{equation}\label{EshkArz_eq1}
T=T_0+T_1+T_2+K, \end{equation} where operators $T_0,$ $T_1,$
$T_2$ and $K$ are defined by the following formulas:
$$T_0f(x,y)=k_0(x,y)f(x,y),$$
$$T_1f(x,y)=\int_{\Omega_1}k_{1}(x,s,y)f(s,y)ds,$$
$$T_2f(x,y)=\int_{\Omega_2}k_{2}(x,t,y)f(x,t)dt,$$
$$Kf(x,y)= \int_{\Omega_1}\int_{\Omega_2}k(x,y;s,t)f(s,t)dsdt.$$
Here $k_0,$ $k_1,$ $k_2$ and $k$ are given measurable functions on
$\Omega_1\times\Omega_2,$ $\Omega_1^2\times\Omega_2^2,$
$\Omega_1\times\Omega_2$ and $(\Omega_1\times\Omega_2)^2,$
respectively, and all integrals have to be understood in the
Lebesgue sense, where $ds=d\mu_{1}(s),dt=d\mu_{2}(t),$
$\mu_{k}(\cdot)$ is the Lebesgue measure on the $\sigma$-algebra
of subsets $\Omega_{k},$ $k=1,2.$

In 1975, Likhtarnikov and Vitova \cite{likh1} spectral properties
of partial integral operators are studied . In \cite{likh1}, the
following restrictions were imposed: $k_1(x,s)\in
L_2(\Omega_1\times\Omega_1),$ $k_2(y,t)\in
L_2(\Omega_2\times\Omega_2)$ and $T_0=K=0.$ In \cite{kal3},
spectral properties of PIO with positive kernels were studied
(under restriction $T_0=K=0$). In Kalitvin and Zabrejko
\cite{kal4}  spectral properties of PIO with kernels of two
variables in $L_p, p\geq1$ are studied. In
\cite{esh4,kal5,kal6,kal7,zab1} for the more general equations
with continuous kernels or kernels in $C(L_1)$ spectral properties
of the PIO and solvability of partial integral equations in the
space $C([a,b]\times[c,d])$ were studied.

Self-adjoint PIO with $T_0\neq0$ were first studied in
\cite{esh1}, where theorem above essential spectrum is proved.
Finiteness and infiniteness of a discrete spectrum of self-adjoint
PIO, arising in the theory of Schr\"{o}dinger operators, are
investigated in \cite{lak1,ras1,esh3}. In \cite{app2} applications
of partial integral equations and operators to problems of
continuous mechanics, elasticity problems and other problems were
considered. Still, some important spectral properties of PIO in
the space $L_2$ are left open. The present paper is dedicated to
the mentioned problem for PIO with degenerate kernels from the
class $L_2$.

Let $T_1$ be a linear integral operator in the space
$L_2([a,b]\times[c,d])$ given by the formula
\begin{equation}\label{EshkArzeq_2}
(T_1f)(x,y)=\int_a^bk(x,s,y)f(s,y)ds.
\end{equation}
Here $k(x,s,y)$ is a measurable function on $[a,b]^2\times[c,d].$

The kernel $k(x,s,y)$ of the integral operator $T_1$ usually
satisfies the condition
$$\int^b_ak(x,s,y)f(s,y)ds\in L_2([a,b]\times[c,d]), \quad \forall f\in
L_2([a,b]\times[c,d]).$$ Consequently, the operator $T_1$ is a
linear bounded operator on $L_2([a,b]\times[c,d]).$ If, in
addition, the kernel $k(x,s,y)$ satisfies the condition:
$$k(x,s,y)=\overline{k(s,x,y)},\quad \ \text{for almost all} \
y\in[c,d],$$ then the operator $T_1$ is a self-adjoint operator on
the Hilbert space $L_2([a,b]\times[c,d]).$

Let $\{\varphi_k(x)\}_{k=1}^n$  be a orthonormal system of
functions from the $L_2[a,b],$ and let $\{h_k(y)\}_{k=1}^n$ be a
system of essential bounded real functions on $[c,d].$

We define a measurable function $k_1(x,s,y)$ on
$[a,b]^2\times[c,d]$ by the following rule:
\begin{equation}\label{EshkArz_eq3}
k_1(x,s,y)=\sum\limits_{k=1}^n\varphi_k(x)\overline{\varphi_k(s)}h_k(y).
\end{equation}
Then the PIO $T_1$ with the kernel $k_1(x,s,y)$ is a self-adjoint
bounded linear operator on $L_2([a,b]\times[c,d]).$

Let $\{\psi_k(y)\}_{k=\overline{1,m}}$ be a some orthonormal
system of functions from the $L_2[c,d],$ $\{p_k
(x)\}_{k=\overline{1,m}}$ be a system of essential bounded real
functions on $[a,b].$ We define the measurable function
$k_2(x,t,y)$ on $[a,b]\times[c,d]^2$ by the following rule:
\begin{equation}\label{EshkArz_eq4}
k_2(x,t,y)=\sum\limits_{j=1}^m
p_j(x)\psi_j(y)\overline{\psi_j(t)}. \end{equation}

Then the PIO $T_2$ with the kernel $k_2(x,t,y):$
$$T_2f(x,y)=\int_c^dk_2(x,t,y)f(x,t)d\mu_2(t)$$ is a linear bounded
self-adjoint operator on $L_2([a,b]\times[c,d]).$

For an essential bounded function $\varphi\geq0$ on the measurable
set $\Omega\subset\mathbb{R}^{\nu},$ we define
$$
esssup_{\Omega}(\varphi)=inf\{C:\mu(\{\xi\in
\Omega:\varphi(\xi)>C\})=0\},
$$
where $\mu(\cdot)$ is the Lebesgue measure on $\mathbb{R}.$ For a
measurable function $\varphi$ on the set
$\Omega\subset\mathbb{R}^{\nu},$ the number $\lambda\in\mathbb{R}$
is called an essential value of the function $\varphi$ if
$$
\mu(\{\xi\in\Omega:
\lambda-\varepsilon<\varphi(\xi)<\lambda+\varepsilon\})>0
$$
for all $\varepsilon>0.$ We denote by $Essran(\varphi)$ the set of
all essential values of the function $\varphi.$

In this paper, we study essential and discrete spectra of PIO of
the form $T_1+T_2$ with the degenerate kernels. The resolvent set,
spectrum, essential spectrum and discrete spectrum are denoted by
$\rho,$ $\sigma,$ $\sigma_{ess}$ and $\sigma_{disc},$,
respectively (see \cite{RS}).

\section{Spectral property of PIO $T_1$ and $T_2$}\label{EshkArz_Section2}

In this section, we study spectra of PIO $T_1$ and $T_2.$

\begin{prop}\label{EshArz_prop2.1.} Zero is an eigenvalue of $T_1$ of infinite multiplicity.
A number $\lambda_0\neq0$ is an eigenvalue of $T_1$ if and only if
there exists $1\leq j_0\leq n$ such that
$\mu_{2}(h_{j_0}^{-1}(\{\lambda_0\}))>0.$
\end{prop}
We denote by $M$ the subspace of a Hilbert space $L_2[a,b]$
constructed by the orthogonal system $\{\varphi_1, ... ,
\varphi_n\}.$ Then $\dim M=n,$ and for the subspace
$\mathcal{H}=L_2[a,b]\circleddash M,$ we have
$\dim\mathcal{H}=\infty.$

Let $\{g_k\}_{k\in\mathbb{N}}$ be an orthonormal basis in
$\mathcal{H},$ and $\psi\in L_2[c,d],$ $\|\psi\|=1.$ Then
$f_k(x,y)=g_k(x)\psi(y)\in L_2([a,b]\times[c,d]),$
$k\in\mathbb{N},$ and the system $\{f_k\}_{k\in\mathbb{N}}$ is
orthonormal in $L_2([a,b]\times[c,d]).$ Clearly,
$$T_1f_k=\sum\limits_{i=1}^n\int_a^b\varphi_i(x)\overline{\varphi_i(s)}f_k(s,y)d\mu_1(s)=0, \,\
k\in\mathbb{N},$$ i.e. zero is an eigenvalue of PIO $T_1$ of
infinite multiplicity.\\
``if''  part. Let $\lambda_0\in\mathbb{C}\backslash\{0\}$ be an
eigenvalue of the PIO $T_1.$ Then there exists $f\in
L_2([a,b]\times[c,d]),$ $\|f\|=1$ such that
$$T_1f=\lambda_0f.$$ We define compact self-adjoint integral
operators $K_\omega$ in $L_2[a,b]$ as follows:
$$K_\omega\varphi(x)=\int_a^bk_1(x,s,w)\varphi(s)d\mu_1(s), \,\
\omega\in\Omega_0,$$ where
$$\Omega_0=\{\omega\in[c,d]: p_{\omega}(x,s)=k_1(x,s,\omega)\in
L_2([a,b]^2)\}.$$ We have $\mu_2([c,d]\setminus\Omega_0)=0$. Put
$$\mathcal{M}_{1}=\{\omega\in\Omega_0: f_{\omega}(x)=f(x,\omega)\in
L_2[a,b]\}.$$ Then $\mu_{2}(\mathcal{M}_{1})>0$.

We define measurable subsets: $$\mathcal{D}_k=\{\omega\in[c,d]:
h_k(\omega)=\lambda_0\}, \,\ k\in\{1, ... ,n\}.$$ Put
$$\mathcal{D}_0=\bigcup\limits_{k=1}^n\mathcal{D}_k.$$

Let $\omega\in \mathcal{M}_{1}.$ Then
$K_{\omega}f_{\omega}=\lambda_0f_{\omega},$  i.e. the number
$\lambda_0$ is an eigenvalue of operators $K_{\omega}, \,\
\omega\in \mathcal{M}_{1}.$ We define by
$\left\{\alpha_1^{(\omega)}, ..., \alpha_{n_{\omega}}^{(\omega)}
\right\}$ the set of eigenvalues of the operator $K_{\omega}$
which are different from zero. Then
$$\lambda_0\in\ \{\alpha_1^{(\omega)}, ...
,\alpha_{n_{\omega}}^{(\omega)}\}\subset\{h_1(\omega), ... ,
h_n(\omega)\},$$  i.e. there exists $j_0\in\{1, ... , n\}$ such
that $h_{j_{0}}(\omega)=\lambda_0.$

Consequently, we have $\omega\in\mathcal{D}_0.$ Thus,
$\mathcal{M}_{1}\subset\mathcal{D}_0.$ It means that
$\mu_2(\mathcal{D}_0)>0.$ Then there exists $j_0\in\{1, ... , n\}$
such that
$$\mu_2(\{\omega\in[c,d]: h_{j_0}(\omega)=\lambda_0\})>0.$$

\noindent ``only if''  part. Let for a $j_{0}\in\{1, ... , n\}$ we
have $\mu_2(h_{j_{0}}^{-1}(\{\lambda_0\}))>0.$ Put
$\mathcal{D}=h_{j_{0}}^{-1}(\{\lambda_0\}).$ We define the
function $\psi\in L_2[c,d]$ in the following way:
$$\psi(y)=\frac{\chi_{\mathcal{D}}(y)}{\sqrt{\mu_2(\mathcal{D})}}, \,\
y\in[c,d],$$ where $\chi_{G}(\cdot)$ is the characteristic
function of a set $G.$ Obviously, $\|\psi\|=1.$ Let
$f(x,y)=\varphi_{j_0}(x)\psi(y).$ Then $f\in
L_2([a,b]\times[c,d])$ and $\|f\|=1.$ On the other hand,
$$T_1f(x,y)=\sum\limits_{k=1}^n\varphi_k(x)\int_a^b\overline{\varphi_k(s)}h_k(y)
\varphi_{j_{0}}(s)\psi(y)d\mu_1(s)=\varphi_{j_{0}}(x)h_{j_{0}}(y)\psi(y)=\lambda_0f(x,y),$$
i.e. the number $\lambda_0$ is an eigenvalue of $T_1.$

We consider the following projectors $P_k$ in the space
$L_2([a,b]\times[c,d]):$
$$P_kf(x,y)=\int_a^b\varphi_{k}(x)\overline{\varphi_k(s)}f(s,y)d\mu_1(s), \,\ k\in\{1, ... , n\}.$$

Let $P=P_1+ ... +P_n$ and $P_0=E-P,$ where $E$ is the identical
operator. Then projectors $P_i$ and $P_j$ $(i\neq j)$ are
orthogonal.

\begin{prop}\label{EshkArz_prop2.2.}
If $\lambda\neq 0$ and $\lambda\overline{\in}
\bigcup\limits_{k=1}^nEssran(h_k),$ then the operator $T_1-\lambda
E$ is invertible in $L_2([a,b]\times[c,d]),$ and the operator
$(T_1-\lambda E)^{-1}$ is bounded in $L_2([a,b]\times[c,d]),$
moreover $$(T_1-\lambda
E)^{-1}f(x,y)=-\frac{1}{\lambda}\left(f(x,y)-\sum\limits_{k=1}^n\frac{h_k(y)}{h_k(y)-\lambda}P_kf(x,y)\right).$$
\end{prop}
\emph{Proof.} Let $\lambda\neq 0$ and
$\lambda\overline{\in}\bigcup\limits_{k=1}^nEssran(h_k).$ We
define the operator $B_{\lambda}$ on the $L_2([a,b]\times[c,d])$
by the formula:
$$B_{\lambda}f(x,y)=\sum\limits_{k=1}^n\frac{1}{h_k(y)-\lambda}P_kf(x,y)-\frac{1}{\lambda}P_0f(x,y).$$
It is clear,
\begin{equation}\label{EshkArz_eq5}
(T_1-\lambda E)B_{\lambda}=B_{\lambda}(T_1-\lambda
E)=E.\end{equation} For all $\lambda\overline{\in}Essran(h_k)
\cup\{0\},$ the operator
$$A_kf(x,y)=\frac{1}{h_k(y)-\lambda}P_kf(x,y), \,\
f\in L_2([a,b]\times[c,d])$$ is a bounded operator in
$L_2([a,b]\times[c,d]).$ Then the operator $B_{\lambda}$ is a
bounded operator in $L_2([a,b]\times[c,d]).$ For each
$\lambda\bar{\in}\{0\}\cup\left(\bigcup\limits_{k=1}^nEssran(h_k)\right),$
\eqref{EshkArz_eq5} implies $$(T-\lambda E)^{-1}=B_\lambda.$$
Hence, we have
$$B_\lambda f(x,y)=-\frac{1}{\lambda}\left(f(x,y)-\sum\limits_{k=1}^n\frac{h_k(y)}{h_k(y)-\lambda}P_kf(x,y)\right).$$

\begin{thm}\label{EshkArz_thm2.3.}
For the spectra $\sigma(T_1)$ of the PIO $T_1$ with a degenerate
kernel \eqref{EshkArz_eq3}, the following formula holds:
$$\sigma(T_1)=\{0\}\cup\left(\bigcup\limits_{k=1}^nEssran(h_k)\right).$$
\end{thm}

\emph{Proof.} By proposition \ref{EshkArz_prop2.2.}, we obtaine
$$\sigma(T_1)\subset\{0\}\cup\left(\bigcup\limits_{k=1}^nEssran(h_k)\right).$$
However, by proposition \ref{EshArz_prop2.1.} we have
$0\in\sigma(T_1).$ Now we prove
$$\bigcup\limits_{k=1}^nEssran(h_k)\subset\sigma(T_1).$$

Let $\lambda_0\in Essran(h_{j_{0}}),$ $\lambda_0\neq0$ and $t_0$
be arbitrary point from the subset
$h_{j_{0}}^{-1}(\{\lambda_0\}).$

Put $$V_i=\left\{t\in[c,d]:
\frac{1}{i+1}<|t_0-t|<\frac{1}{i}\right\}, \,\ i\in\mathbb{N}.$$ \\
Then there exists $n_0\in\mathbb{N}$ such that $\mu_2(V_i)>0$ for
all $i\geq n_0.$ We consider the following sequence of orthonormal
functions $\chi_p(y)\in L_2[c,d]:$

\[ \chi_p(y) = \left\{
\begin{array}{ll}
\frac{1}{\sqrt{\mu_2(V_p)}}, & y\in V_p, \\ \\

0, & y\overline{\in} V_p,
\end{array} \right. \] \\
where $p\geq n_0.$ We define by $f_p(x,y)\in L_2([a,b]\times
[c,d])$ the orthonormal system of functions:
$f_p(x,y)=\varphi_{j_{0}}(x)\chi_p(y), \,\ p\geq n_0.$ Then we
have
$$(T_1-\lambda_0E)f_p(x,y)=(h_{j_{0}}(y)-\lambda_0)f_p(x,y).$$
Hence,
$$\|(T_1-\lambda_0E)f_p\|\leq\sqrt{esssup_{V_p}\left(h_{j_{0}}(y)-\lambda_0\right)^2}, \,\ p\geq
n_0.$$ Since zero is an essential value of the function
$h_{j_{0}}(y)-\lambda_0,$ then for large $n_1\geq n_0,$ there
exists a small number $\delta_{n_{1}}$ such that
$$\mid h_{j_{0}}(y)-\lambda_0\mid<\delta_{n_{1}} \,\,\ \mbox{for almost all} \,\,\ y\in V_p, \,\ p\geq
n_1.$$ Therefore
$$\parallel(T_1-\lambda_0E)f_p\parallel<\delta_{n_{1}}, \,\,\ p\geq
n_1,$$ i.e.
$\lim\limits_{n\rightarrow\infty}\parallel(T_1-\lambda_0E)f_n\parallel=0.$
This and the Weyl criterion for an essential spectrum of
self-adjoint operators \cite{RS} imply
$\lambda_0\in\sigma_{ess}(T_1)\subset\sigma(T_1).$

\begin{prop}\label{EshkArz_prop2.4.}
Any eigenvalue of the PIO $T_1$ is of infinite multiplicity.
\end{prop}
\emph{Proof.} Let $\lambda\in\mathbb{R}\backslash\{0\}$ be an
eigenvalue of $T_1.$ Then there exists $f_0\in
L_2([a,b]\times[c,d]), \,\ \|f_0\|=1$ such that $T_1f_0=\lambda
f_0.$ We define the measurable subset $\Omega_0\subset[c,d]:$
$$\Omega_0=\left\{y\in[c,d]: \int_a^b|f_0(x,y)|^2d\mu_1(x)\neq0\right\}.$$
Obviously, $\mu_2(\Omega_0)>0.$ Define the function:

\[ f_0(x,y)= \left\{
\begin{array}{ll}
\frac{f_0(x,y)}{\sqrt{\int_a^b|f_0(s,y)|^2d\mu_1(s)}}, & x\in [a,b], y\in\Omega_0 \\ \\

 0, & x\in[a,b], y\overline{\in}\Omega_0,
\end{array} \right. \]
Then $f_0\in L_2([a,b]\times[c,d]),$ and $f_0\neq 0.$

Let $\{\widetilde{\psi_k}\}_k\in\mathbb{N}$ be a system of
orthonormal functions from $L_2(\Omega_0).$ Consider the sequence
of functions from $L_2([a,b]\times[c,d]):$
$$f_k(x,y)=f_0(x,y)\psi_k(y), \,\ k\in\mathbb{N},$$ where

\[\psi_k(y) = \left\{
\begin{array}{ll}
\widetilde{\psi_k}(y), & y\in\Omega_0 \\ \\
0, & y\overline{\in}\Omega_0.
\end{array} \right. \]

Then
$$\int_a^b\int_c^d\left|f_n(x,y)\right|^2d\mu_1(x)d\mu_2(y)=\int_{\Omega_{0}}\left|\widetilde{\psi_n}(y)\right|^2d\mu_2(y)=1, $$
and
$$\left(f_i,f_j
\right)=\int_{\Omega_{0}}\widetilde{\psi_i}(y)\overline{\widetilde{\psi_j}(y)}d\mu_2(y)=0$$
for $i\neq j.$

Clearly $$T_1f_k(x,y)=\lambda f_k(x,y), \,\ k\in\mathbb{N},$$ i.e.
the number $\lambda$ is an eigenvalue of the PIO $T_1$ of infinite
multiplicity.

\begin{cor}\label{EshArz_cor2.5.} A discrete spectrum of the PIO $T_1$ with a
degenerate kernel \eqref{EshkArz_eq3} is absent.  \end{cor}

\begin{cor}\label{EshArz_cor2.6.} If every function $h_k,$ $k\in\{1, ...
, n\}$ is continuous and strictly monotone on $[c,d],$ then there
is not an eigenvalue of the PIO $T_1$ different from
zero.\end{cor} Consider the following projectors $Q_j$ in the
space $L_2([a,b]\times[c,d]):$
$$Q_jf(x,y)=\int_c^d\psi_j(y)\overline{\psi_j(t)}f(x,t)d\mu_2(t), \,\,\ j\in\{1, ... , m\}.$$

\begin{prop}\label{EshArz_prop2.7.} If $\lambda\neq0$ and
$\lambda\overline{\in}\bigcup\limits_{j=1}^m Essran(p_j),$ then
the operator $T_2-\lambda E$ is invertible on
$L_2([a,b]\times[c,d]),$ and the operator $(T_2-\lambda E)^{-1}$
is bounded in $L_2([a,b]\times[c,d]),$ moreover
$$\left(T_2-\lambda E\right)^{-1}f(x,y)=-\frac{1}{\lambda}\left(f(x,y)-\sum\limits_{j=1}^m\frac{p_j(x)}{p_j(x)-\lambda}Q_jf(x,y)\right).$$
\end{prop}

\begin{thm}\label{EshkArz_thm2.8.}
For the spectrum $\sigma(T_2)$ of the PIO $T_2$ with a degenerate
kernel \eqref{EshkArz_eq4}, the following formula is hold:
$$\sigma(T_2)=\{0\}\cup\left(\bigcup\limits_{j=1}^m Essran(p_j)\right).$$ \end{thm}

\section{Solvability of partial integral equations}\label{EshkArz_Section3}

We consider the Fredholm partial integral equation (PIE) of the
second kind
\begin{equation}\label{EshkArz_eq6}
f(x,y)-\tau(T_1+T_2)f(x,y)=g(x,y) \end{equation} in the Hilbert
space $L_2([a,b]\times[c,d]),$ where $f$ is an unknown function
from $L_2([a,b]\times[c,d]),$ $g\in L_2([a,b]\times[c,d])$ is a
given function, and $\tau\in\mathbb{C}$ is the parameter of the
equation.

The homogeneous PIE corresponding to \eqref{EshkArz_eq6} has the
following form:
$$f(x,y)-\tau(T_1+T_2)f(x,y)=0.$$

In this section, we reduce to some necessary results for PIE of
the second kind.

Assume that $\tau\neq0$ and $\tau^{-1}\in\rho(T_1).$ Then the
operator $E-\tau T_1$ is invertible in $L_2([a,b]\times[c,d]),$
and the operator $(E-\tau T_1)^{-1}$ is bounded on
$L_2([a,b]\times[c,d]),$ moreover, by proposition
\ref{EshkArz_prop2.2.}, we have
$$(E-\tau T_1)^{-1}f(x,y)=f(x,y)+\tau\sum\limits_{k=1}^n\frac{h_k(y)}{1-\tau h_k(y)}P_kf(x,y).$$

Analogously, if $\tau\neq0$ and $\tau^{-1}\in\rho(T_2),$ then the
operator $E-\tau T_2$ is invertible in $L_2([a,b]\times[c,d]),$
and the operator $(E-\tau T_2)^{-1}$ is bounded on
$L_2([a,b]\times[c,d]),$ moreover, by proposition
\ref{EshArz_prop2.7.},
$$(E-\tau T_2)^{-1}f(x,y)=f(x,y)+\tau\sum\limits_{j=1}^m\frac{p_j(x)}{1-\tau
p_j(x)}Q_jf(x,y).$$

Let $\tau\neq0$ and $\tau^{-1}\in\rho(T_1)\cap\rho(T_2).$ We
define compact operators $W_1(\tau)$ and $W_2(\tau)$ by the
following formulas:
$$W_1(\tau)=(E-\tau T_2)^{-1}S_1(\tau)T_2, \,\,\ W_2(\tau)=(E-\tau T_1)^{-1}S_2(\tau)T_1,$$
where
$$S_1(\tau)f(x,y)=\sum\limits_{k=1}^n\frac{h_k(y)}{1-\tau
h_k(y)}(P_kf)(x,y), \,\
S_2(\tau)f(x,y)=\sum\limits_{i=1}^m\frac{p_i(x)}{1-\tau
p_i(x)}(Q_if)(x,y).$$

\begin{lem}\label{EshArz_lem3.3.} Let $\tau\neq0$ and
$\tau^{-1}\in\rho(T_1)\cap\rho(T_2).$ Then the following three
homogenous Fredholm PIE of the second kind are equivalent:
\begin{equation}\label{EshkArz_eq7}
f-\tau(T_1+T_2)f=0,  \end{equation}
\begin{equation}\label{EshkArz_eq8}
f-\tau^2W_1(\tau)f=0,  \end{equation}
\begin{equation}\label{EshkArz_eq9}
f-\tau^2W_2(\tau)f=0.  \end{equation}
\end{lem}

In \cite{esh4,kal5,kal6,kal7, zab1} lemmas similar to the lemma
\ref{EshArz_lem3.3.} for the case of general PIE's in
$C([a,b]\times[c,d])$ with continuous kernels or kernels in
$C(L_1)$ were proved. The scheme for the proof of lemma
\ref{EshArz_lem3.3.} can be seen from these works.

Assume that $\tau\neq0$ and $\tau^{-1}\in\rho(T_1)\cap\rho(T_2).$
We denote by $\Delta_1(\tau)$ and $\Delta_2(\tau)$ the Fredholm
determinants of the operators $E-\tau^{2}W_1(\tau)$ and
$E-\tau^{2}W_2(\tau),$ respectively. Define in $\mathbb{C}$ the
following subsets
$$\mathcal{R}_1=\left\{\tau\in\mathbb{C}\setminus\{0\}: \tau^{-1}\in\rho(T_1)\cap\rho(T_2) \,\,\  \mbox{and} \,\,\ \Delta_1(\tau)\neq0\right\},$$
$$\mathcal{R}_2=\left\{\tau\in\mathbb{C}\setminus\{0\}: \tau^{-1}\in\rho(T_1)\cap\rho(T_2) \,\,\  \mbox{and} \,\,\ \Delta_2(\tau)\neq0\right\},$$
$$\mathcal{D}_1=\left\{\tau\in\mathbb{C}\setminus\{0\}: \tau^{-1}\in\rho(T_1)\cap\rho(T_2) \,\,\  \mbox{and} \,\,\ \Delta_1(\tau)=0\right\},$$
$$\mathcal{D}_2=\left\{\tau\in\mathbb{C}\setminus\{0\}: \tau^{-1}\in\rho(T_1)\cap\rho(T_2) \,\,\  \mbox{and} \,\,\ \Delta_2(\tau)=0\right\}.$$

It follows from lemma \ref{EshArz_lem3.3.} that
$\mathcal{R}_1=\mathcal{R}_2$ and $\mathcal{D}_1=\mathcal{D}_2.$
Put
$$\mathcal{R}=\mathcal{R}(T)=\mathcal{R}_1 \,\,\ \mbox{and} \,\,\ \mathcal{D}=\mathcal{D}(T)=\mathcal{D}_1.$$
Then we obtain
\begin{equation}\label{EshkArz_eq10}\mathcal{R}\cup\mathcal{D}=\left\{\tau\in\mathbb{C}: \tau\in\mathbb{C}\setminus\{0\}  \,\,\ \mbox{and} \,\,\
\tau^{-1}\in\rho(T_1)\cap\rho(T_2)\right\}\end{equation} and
$\mathcal{R}\cap\mathcal{D}=\emptyset.$

Lemma \ref{EshArz_lem3.3.} implies the following
\begin{thm}\label{EshkArz_thm3.4.}
Let $\tau\in\mathcal{R}\cup\mathcal{D}.$ Homogenous PIE
\eqref{EshkArz_eq7} has a non-trivial so\-lu\-tion if and only if
$\tau\in\mathcal{D}.$ \end{thm}

\begin{thm}\label{EshkArz_thm3.5.}
Let $\tau\in\mathcal{R}.$ Then the Fredholm PIE of the second kind
\begin{equation}\label{EshkArz_eq11} f-\tau(T_1+T_2)f=g
\end{equation} has the unique solution $f_0\in
L_2([a,b]\times[c,d])$  for any $g\in L_2([a,b]\times[c,d]).$
\end{thm}

\section{Spectrum of the PIO $T_1+T_2$}\label{EshkArz_Section4}

We put
$$\mathcal{D}_0=\mathcal{D}_0(T)=\left\{\xi:  \xi=\frac{1}{\tau}, \,\
\tau\in\mathcal{D}(T)\right\}.$$\\

\begin{lem}\label{EshArz_lem4.1.}For each
$\lambda\in\mathbb{C} \setminus \left(\sigma(T_1)\cup
\sigma(T_2)\cup\mathcal{D}_0(T)\right),$ the resolvent
$R_{\lambda}(T)$ of the PIO $T=T_1+T_2$ exists and is bounded on
$L_2([a,b]\times[c,d]).$ \end{lem}

\emph{Proof.} Let $\lambda\in\mathbb{C}\setminus\left( \sigma(T_1)
\cup \sigma(T_2)\cup\mathcal{D}_0(T)\right).$ Then $\lambda\neq0,$
and $\lambda\in\rho(T_1)\cap\rho(T_2).$ Then operators
$E-\frac{1}{\lambda}T_1$ and $E-\frac{1}{\lambda}T_2$ are
injective and
$$\left(E-\frac{1}{\lambda}T_1\right)^{-1}=E+\frac{1}{\lambda}S_1\left(\frac{1}{\lambda}\right),\,\,\
\left(E-\frac{1}{\lambda}T_2\right)^{-1}=E+\frac{1}{\lambda}S_2\left(\frac{1}{\lambda}\right).$$ \\
However, from $\lambda\overline{\in}\mathcal{D}_0(T)$ and by
\eqref{EshkArz_eq10}, we obtain
$\Delta_1\left(\frac{1}{\lambda}\right)\neq0.$ Consequently, the
operator $E-\frac{1}{\lambda^2}W_1\left(\frac{1}{\lambda}\right)$
is injective. By the other hand, we have
$$T_1+T_2-\lambda E=(T_1-\lambda E)(E+(T_1-\lambda E)^{-1}T_2)=$$
$$=(T_1-\lambda
E)\left[E-\frac{1}{\lambda}\left(E+\frac{1}{\lambda}S_1\left(\frac{1}{\lambda}\right)\right)T_2\right]=$$
$$=(T_1-\lambda
E)\left(E-\frac{1}{\lambda}T_2\right)\left(E-\frac{1}{\lambda^2}\left(E-\frac{1}{\lambda}T_2\right)^{-1}S_1\left(\frac{1}{\lambda}\right)T_2\right)=$$
$$=-\lambda\left(E-\frac{1}{\lambda}T_1\right)\left(E-\frac{1}{\lambda}T_2\right)\left(E-\frac{1}{\lambda^2}W_1\left(\frac{1}{\lambda}\right)\right).$$

Hence,
$$R_{\lambda}(T)=-\frac{1}{\lambda}\left(E-\frac{1}{\lambda^2}W_1\left(\frac{1}{\lambda}\right)\right)^{-1}\left(E-\frac{1}
{\lambda}T_2\right)^{-1}\left(E-\frac{1}{\lambda}T_1\right)^{-1}.$$
Boundedness of the operator $R_{\lambda}(T)$ follows from the last
equality.

\begin{cor}\label{EshArz_cor4.2.} For the resolvent operator $R_{\lambda}(T)$ of the PIO
$T=T_1+T_2,$ the formula
$$R_{\lambda}(T)=-\frac{1}{\lambda}\left(E-\frac{1}{\lambda^2}W_1\left(\frac{1}{\lambda}\right)\right)^{-1}\left(E-\frac{1}
{\lambda}T_2\right)^{-1}\left(E-\frac{1}{\lambda}T_1\right)^{-1}$$
holds for each $\lambda\in\mathbb{C}\setminus\left(\sigma(T_1)
\cup\sigma(T_2)\cup\mathcal{D}_0 (T)\right)$.  \end{cor}

\begin{prop}\label{EshArz_prop4.3.} Zero is an eigenvalue of infinite multiplicity of the PIO $T=T_1+T_2.$
\end{prop}

\emph{Proof.} Define by $\mathcal{L}$ the subspace of a Hilbert
space $L_2([a,b]\times[c,d])$ constructed by the orthogonal system
$\left\{\varphi_i(x)\psi_j(y)\right\}_{i=\overline{1,n},\
j=\overline{1,m}}$. Then $\dim\mathcal{L}=m\times n$ and for
subspace $\mathcal{H}_0=L_2([a,b]\times[c,d])\ominus\mathcal{L},$
we have $\dim\mathcal{H}_0=\infty.$

Let $\left\{f_k\right\}_{k\in\mathbb{N}}$ be an orthonormal basis
in $\mathcal{H}_0$. It is obvious that
$$T_1f_p(x,y)=\sum\limits_{i=1}^n\int_a^b\varphi_i(x)\overline{\varphi_i(s)}h_i(y)f_k(s,y)d\mu_1(s)=0, \,\ p\in\mathbb{N},$$
$$T_2f_p(x,y)=\sum\limits_{j=1}^m\int_c^dp_j(x)\psi_j(y)\overline{\psi_j(t)}f_k(x,t)d\mu_2(t)=0, \,\ p\in\mathbb{N}.$$
i.e. $$(T_1+T_2)f_k(x,y)=0, \,\ k\in\mathbb{N}.$$

\begin{prop}\label{EshArz_prop4.4.} The inclusion
$\sigma(T_1)\cup\sigma(T_2)\subset\sigma_{ess}(T_1+T_2)$
holds.\end{prop}

\emph{Proof.} We show that
$\sigma(T_1)\subset\sigma_{ess}(T_1+T_2)$ (the inclusion
$\sigma(T_2)\subset\sigma_{ess}(T_1+T_2)$ can be proved
analogously). We have
$$\sigma(T_1)=\{0\}\cup\left(\bigcup\limits_{i=1}^n
Essran(h_i)\right).$$ By proposition \ref{EshArz_prop4.3.},
$0\in\sigma_{ess}(T_1+T_2).$ Assume that $\lambda_0\in
Essran(h_{i_{0}}),$ $\lambda_0\neq0,$ where $i_0\in\{1,2, ... ,
n\}$ and $y_0$ is arbitrary point of the set
$h_{i_{0}}^{-1}\left(\{\lambda_0\}\right).$ Put
$$V_k=\left\{y\in[c,d]: \frac{1}{k+1}<|y_0-y|<\frac{1}{k}\right\}, \,\ k\in\mathbb{N}.$$
Then there exists $k_0\in\mathbb{N}$ such that $\mu_2(V_k)>0$ for
all $k\geq k_0,$ and $\lim\limits_{n\rightarrow\infty}
\mu_2(V_k)=0.$ Consider the sequence of orthonormal functions
$\Phi_k(y)\in L_2[c,d]:$
$$\Phi_k(y)=\frac{\chi_{V_k}(y)}{\sqrt{\mu_2(V_k)}}, \,\,\ k\in\mathbb{N}, \,\ k\geq k_0.$$
Put $$f_k(x,y)=\varphi_{i_0}(x)\Phi_k(y), \,\ k\geq k_0.$$ Then
the system of functions $\left\{f_k(x,y)\right\}_{k\geq k_0}$ of
$L_2([a,b]\times[c,d])$ is an orthonormal system. \\ Now we prove
that
$\lim\limits_{k\rightarrow\infty}\|(T_1+T_2-\lambda_0E)f_k\|=0.$
However, $\lim\limits_{k\rightarrow\infty}
\|(T_1-\lambda_0E)f_k\|=0$ (see the proof of theorem
\ref{EshkArz_thm2.3.}). We show that $\lim\limits_{k\rightarrow
\infty} \|T_2f_k\|=0.$

We define operators $A_i,$ $i=1, ... ,m$ in the following way:
$$A_if(x,y)=\int_c^dp_i(x)\psi_i(y)\overline{\psi_i(t)}f(x,t)d\mu_2(t), \,\ f\in L_2([a,b]\times[c,d]).$$
Then
$$\|A_if_k\|^2\leq\int\limits_a^b\int\limits_c^d\left(\int\limits_c^d\mid
p_i(x)\mid\cdot\mid \psi_i(t)\mid\cdot\mid \psi_i(y)\mid\cdot\mid
\varphi_{i_{0}}(x) \mid\cdot\mid \Phi_k(t)\mid
d\mu_2(t)\right)^2d\mu_1(x)\cdot$$
$$\cdot d\mu_2(y)\leq C_i^2\left(\int\limits_c^d\mid \psi_i(t)\mid\cdot\Phi_k(t)d\mu_2(t)\right)^2=
\frac{C_i^2}{\mu_2(V_k)}\left(\int\limits_{V_{k}}\mid\psi_i(t)\mid
d\mu_2(t)\right)^2\leq$$
$$\leq C_i^2\cdot\int\limits_{V_k}\mid \psi_i(t)\mid^2d\mu(t), \,\,\ i\in\{1, ... ,
n\},$$  where $C_i=esssup_{[a,b]}\mid p_i(x)\mid.$ Since Lebesgue
integrals are absolute continuous, we obtain
$$\lim\limits_{k\rightarrow\infty}\int_{V_k}\mid
\psi_i(t)\mid^2d\mu_2(t)=0$$ at
$\lim\limits_{k\rightarrow\infty}\mu_2(V_k)=0$ and $\|\psi_i\|=1.$

Thus, we get $\lim\limits_{k\rightarrow\infty}\|A_if_k\|=0,$
$i\in\{1, ... , m\}.$

However,
$$\|T_2f_k\|\leq\sum\limits_{i=1}^m\|A_if_k\|,$$
what follows $\lim\limits_{k\rightarrow\infty}\|T_2f_k\|=0.$
Hence, $\lim\limits_{k\rightarrow\infty}\|(T_1+T_2-\lambda_0E)
f_k\|=0.$ Finally, by the Weyl criterion for an essential spectrum
of self-adjoint operators \cite{RS},
$\lambda_0\in\sigma_{ess}(T_1+T_2).$

\begin{prop}\label{EshArz_prop4.5.} Each $\lambda\in\mathcal{D}_0(T)$
is an eigenvalue of finite multiplicity of the PIO $T=T_1+T_2.$
\end{prop}

\emph{Proof.} Let $\lambda\in\mathcal{D}_0(T).$ Then
$\lambda\neq0$ and $\Delta_1\left(\frac{1}{\lambda}\right)=0,$
where $\Delta_1(\tau)$ is the Fredholm determinant of the operator
$E-\tau^2W_1(\tau).$ It means that the number 1 is an eigenvalue
of the compact integral operator
$\frac{1}{\lambda^2}W_1\left(\frac{1}{\lambda}\right).$ By lemma
\ref{EshArz_lem3.3.}, the number $\lambda$ is an eigenvalue of the
PIO $T_1+T_2$. Since the following integral equations
$$f-\frac{1}{\lambda}(T_1+T_2)f=0$$ and
$$f-\frac{1}{\lambda^2}W_1\left(\frac{1}{\lambda}\right)f=0$$ are
equivalent, the number $\lambda$ is an eigenvalue of finite
multiplicity of $T_1+T_2$ because of every eigenvalue
$\alpha\neq0$ of compact operators is of finite multiplicity.

The next theorem follows from lemma \ref{EshArz_lem4.1.} and
propositions \ref{EshArz_prop4.4.}, \ref{EshArz_prop4.5.}.

\begin{thm}\label{EshkArz_thm4.6.}
For the spectrum $\sigma(T)$ of the PIO $T=T_1+T_2$ with a
degenerate kernels, the following formula
$$\sigma(T_1+T_2)=\{0\}\cup\left(\bigcup\limits_{i=1}^n Essran(h_i)\right)
\cup\left(\bigcup\limits_{j=1}^m
Essran(p_j)\right)\cup\mathcal{D}_0(T)$$ holds. \end{thm}

\section{Discrete spectrum of the PIO $T_1+T_2$}\label{EshkArz_Section5}

Put $G=\mathbb{C}\setminus(\sigma(T_1)\cup\sigma(T_2)).$ It is
well-known that spectra of a linear bounded self-adjoint operators
are compact set in the set of all real numbers. Consequently, the
set $\sigma(T_1)\cup\sigma(T_2)$ is a compact subset in
$\mathbb{R}$. Therefore the set $G$ is an open subset in
$\mathbb{C}$ and $G$ is unbounded domain in $\mathbb{C}$.

For each $\lambda\in G,$ we consider the kernel of the compact
integral operator $W_1\left(\frac{1}{\lambda}\right)$ given as
follows:
$$W_1\left(\frac{1}{\lambda}\right)=\left(E-\frac{1}{\lambda}T_2\right)^{-1}S_1\left(\frac{1}{\lambda}\right)T_2.$$\\

\begin{prop}\label{EshArz_prop5.1.} For the kernel
$\mathcal{K}(x,y;s,t|\lambda)$ ($\lambda\in G$) of the Fredholm
integral operator $W_1\left(\frac{1}{\lambda}\right),$ the
equality
\begin{equation}\label{EshkArz_eq13}
\mathcal{K}(x,y;s,t|
\lambda)=\lambda\sum\limits_{j=1}^n\sum\limits_{k=1}^m
F_{k,j}(x,y;\lambda)B_{k,j}(s,t), \end{equation} is valid, where
\begin{equation}\label{EshkArz_eq14}
F_{k,j}(x,y;\lambda)=\varphi_j(x)\!\left(\frac{\psi_k(y)h_j(y)}{\lambda-h_j(y)}+\sum\limits_{i=1}^m\frac{p_i(x)\psi_i(y)}{\lambda-p_i(x)}
\int\limits_c^d\frac{h_j(\xi)}{\lambda-h_j(\xi)}\psi_k(\xi)\overline{\psi_i(\xi)}d\mu_2(\xi)\right),\end{equation}
$$B_{k,j}(s,t)=p_k(s)\overline{\varphi_j(s)}\cdot\overline{\psi_k(t).}$$ \end{prop} \

\emph{Proof.} Let $\lambda\in G.$ Then
$\lambda\in\rho(T_1)\cap\rho(T_2),$ and we get
$$\left(E-\frac{1}{\lambda}T_2\right)^{-1}=E+\frac{1}{\lambda}S_2\left(\frac{1}{\lambda}\right).$$
However,
$$W_1\left(\frac{1}{\lambda}\right)=S_1\left(\frac{1}{\lambda}\right)T_2+\frac{1}{\lambda}
S_2\left(\frac{1}{\lambda}\right)S_1\left(\frac{1}{\lambda}\right)T_2.$$

For each $f\in L_2([a,b]\times[c,d]),$ using representations of
operators $S_1(\tau)$ and $S_2(\tau),$ we obtain
$$S_1\left(\frac{1}{\lambda}\right)T_2f(x,y)=\sum\limits_{j=1}^n\sum\limits_{k=1}^m\int\limits_a^b\int\limits_c^d K_{k,j}(x,y;\lambda)
B_{k,j}(s,t)f(s,t)d\mu_1(s)d\mu_2(t),$$
$$S_2\left(\frac{1}{\lambda}\right)S_1\left(\frac{1}{\lambda}\right)T_2f(x,y)=\sum\limits_{j=1}^n\sum\limits_{k=1}^m
\int\limits_a^b\int\limits_c^d G_{k,j}(x,y;\lambda)
B_{k,j}(s,t)f(s,t)d\mu_1(s)d\mu_2(t),$$ where
$$K_{k,j}(x,y;\lambda)=\frac{\lambda\varphi_j(x)\psi_k(y)h_j(y)}{\lambda-h_j(y)},$$
$$G_{k,j}(x,y;\lambda)=\lambda^2\sum\limits_{i=1}^m\frac{p_i(x)\psi_i(y)}
{\lambda-p_i(x)}\int_c^d\frac{h_j(\xi)}{\lambda-h_j(\xi)}\psi_k(\xi)\overline{\psi_i(\xi)}d\mu_2(\xi).$$

Hence, we obtain equality \eqref{EshkArz_eq13} for the kernel
$\mathcal{K}(x,y;s,t|\lambda)$ of the integral operator
$W_1\left(\frac{1}{\lambda}\right).$

Set $$\Gamma_1=\{1,2, ... ,m\}, \,\,\ \Gamma_2=\{1,2, ... ,n\}
\,\,\ \mbox{and} \,\,\ \Gamma=\Gamma_1\times\Gamma_2.$$

We can introduce the relation of partial order in the set $\Gamma$
by the following way: for elements $\omega=(k_1,j_1)\in\Gamma$ and
$\omega'=(k_2,j_2)\in\Gamma,$ we write $\omega\leq\omega'$ if
$k_1<k_2$ or $k_1=k_2,$ $j_1\leq j_2.$ As the set $\Gamma$ is
finite, the set $\Gamma$ is linear complete ordered, i.e. for
arbitrary $\omega,\omega'\in\Gamma,$ we have $\omega\leq\omega'$
or $\omega'\leq\omega.$ Thus, we can give elements of $\Gamma$ in
the increase order:
$$\Gamma=\{\omega_1, \omega_2, ... , \omega_{m\cdot (n-1)}, \omega_{m\cdot n}\},$$
moreover
$$\omega_1=(1,1)<\omega_2=(1,2)< ...<\omega_n=(1,n)<\omega_{n+1}=(2,1)< ... <\omega_{m\cdot n}=(m,n).$$ \

Let $\lambda\in G$ be fixed. For every $\omega=(k,j)\in\Gamma,$ we
define the function $F_{\omega}(x,y;\lambda)$ on
$[a,b]\times[c,d]$ by the following formula:
$$F_{\omega}(x,y;\lambda)=F_{k,j}(x,y;\lambda).$$

Consider the homogenous Fredholm integral equation
\begin{equation}\label{EshkArz_eq15}
f(x,y)-\frac{1}{\lambda^2}W_1\left(\frac{1}{\lambda}\right)f(x,y)=0,
\,\ f\in L_2([a,b]\times[c,d]). \end{equation} Set
$$\int\limits_a^b\int\limits_c^dF_{\omega_{i}}(x,y;\lambda)f(x,y)d\mu(x)d\mu(y)=A_{\omega_{i}}(\lambda), \,\ i\in\{1, ...
,m\cdot n\}.$$ Then the homogenous equation \eqref{EshkArz_eq15}
turns into the equation
$$f(x,y)=\frac{1}{\lambda}\sum\limits_{i=1}^{m\cdot n}A_{\omega_{i}}(\lambda)F_{\omega_{i}}(x,y;\lambda).$$

Let
$$\int\limits_a^b\int\limits_c^dF_{\omega_{i}}(x,y;\lambda)B_{\omega_{l}}(x,y)d\mu_1(x)d\mu_2(y)=\Pi_{i,l}(\lambda), \,\ i,l\in\{1, ...
,m\cdot n\},$$ where $$B_{\omega_{l}}(x,y)=B_{k_{l},j_{l}}(x,y),
\,\,\ \omega_l=(k_l, j_l).$$  Then we obtain a system of
homogenous linear algebraic equations for unknown numbers
$A_{\omega_i}(\lambda):$
\begin{equation}\label{EshkArz_eq16}
A_{\omega_i}(\lambda)-\frac{1}{\lambda}\sum\limits_{l=1}^{m\cdot
n}\Pi_{i,l}(\lambda)A_{\omega_{l}}(\lambda)=0, \,\ i\in\{1, ...
,m\cdot n\}.\end{equation}

\begin{lem}\label{EshArz_lem5.2.}Let $\lambda\in G.$ The homogenous integral
equation \eqref{EshkArz_eq15} has a nontrivial solution if and
only if $\Delta(\lambda)=0,$ where \end{lem}

\[\Delta(\lambda) =
\begin{vmatrix}

\Pi_{1,1}(\lambda)-\lambda & \Pi_{1,2}(\lambda) & \ldots &
\Pi_{1,m\cdot n}(\lambda) \\ \\
 \Pi_{2,1}(\lambda) &
\Pi_{2,2}(\lambda)-\lambda & \ldots & \Pi_{2,m\cdot n}(\lambda) \\ \\
\ldots & \ldots & \ldots & \ldots \\ \\
 \Pi_{m\cdot n,1}(\lambda) &
\Pi_{m\cdot n,2}(\lambda) & \ldots & \Pi_{m\cdot n,m\cdot
n}(\lambda)-\lambda

\end{vmatrix}. \qquad \] \\

\emph{Proof.} Let $\lambda\in G.$ Then equivalence of the Fredholm
integral equation of the second kind \eqref{EshkArz_eq15} and the
system of linear algebraic homogeneous equations
\eqref{EshkArz_eq16} is clear. The determinant
$\widetilde{\Delta}(\lambda)$ of the system of equations
\eqref{EshkArz_eq16} has the following form:

\[\widetilde{\Delta}(\lambda) =
\begin{vmatrix}

1-\frac{\Pi_{1,1}(\lambda)}{\lambda} &
-\frac{\Pi_{1,2}(\lambda)}{\lambda} & \ldots &
-\frac{\Pi_{1,m\cdot n}(\lambda)}{\lambda} \\ \\

 -\frac{\Pi_{2,1}(\lambda)}{\lambda} &
1-\frac{\Pi_{2,2}(\lambda)}{\lambda} & \ldots &
-\frac{\Pi_{2,m\cdot n}(\lambda)}{\lambda} \\ \\

\ldots & \ldots & \ldots & \ldots \\ \\

-\frac{\Pi_{m\cdot n,1}(\lambda)}{\lambda} & -\frac{\Pi_{m\cdot
n,2}(\lambda)}{\lambda} & \ldots &
1-\frac{\Pi_{m\cdot n,m\cdot n}(\lambda)}{\lambda} \\ \\

\end{vmatrix}, \qquad \]
and
$$\widetilde{\Delta}(\lambda)=\left(-\frac{1}{\lambda}\right)^{m\cdot n}\Delta(\lambda).$$

It is well-known, the system of linear homogeneous equations
\eqref{EshkArz_eq16} has nontrivial solution if and only if
$\widetilde{\Delta}(\lambda)=0,$ i.e. $\Delta(\lambda)=0.$
However, we obtain that the homogeneous Fredholm equation
\eqref{EshkArz_eq15} has nontrivial solution if and only if
$\Delta(\lambda)=0.$

\begin{lem}\label{EshArz_lem5.3.} The function $\Delta(z)$ is holomorphic in the
domain $G.$ \end{lem}

\emph{Proof.} Let $\omega\in\Gamma.$ It is known, the function
$F_{\omega}(\lambda)=F_{\omega}(x,y;\lambda)$ is holomorphic by
$\lambda$ in the domain $G$ for almost all $(x,y)\in
[a,b]\times[c,d],$ and for every $\lambda\in G$ the integral
$$\int\limits_a^b\int\limits_c^dF_{\omega}(x,y;\lambda)B_{\omega'}(x,y)d\mu_1(x)d\mu_2(y), \,\,\
\omega,\omega'\in\Gamma$$ exists and is finite. Then for every
$\omega=(i,l)\in\Gamma,$ the function $\Pi_{i,l}(z)$ is a
holomorphic function in $G$. Consequently, the function
$\Delta(z)$ is a sum  of holomorphic functions $F_{\omega_1}(z),
F_{\omega_2}(z), ..., F_{\omega_{m\cdot n}}(z),$ i.e. $\Delta(z)$
is holomorphic in $G.$

\begin{rem}\label{EshArz_rem5.4.} An analogue of Lemma \ref{EshArz_lem5.3.} can be
proved for the general PIE. \end{rem}

\begin{thm}\label{EshkArz_thm5.5.}
The discrete spectrum of the PIO $T=T_1+T_2$ coincides with the
set $\mathcal{D}_0(T).$ \end{thm}

\emph{Proof.} Lemmas \ref{EshArz_lem3.3.} and \ref{EshArz_lem5.2.}
imply
$$\mathcal{D}_0(T)=\{\lambda\in G: \Delta(\lambda)=0\}.$$
By proposition \ref{EshArz_prop4.4.}, we have
$$\sigma(T_1)\cup\sigma(T_2)\subset\sigma_{ess}(T_1+T_2).$$
By theorem \ref{EshkArz_thm4.6.}, we obtain
$$\sigma_{disc}(T)\subset\mathcal{D}_0(T).$$

Let $\lambda_0\in\mathcal{D}_0(T)$ be arbitrary. Then by
proposition \ref{EshArz_prop4.5.}, the number $\lambda_0$ is an
eigenvalue of finite multiplicity of the operator $T_1+T_2.$ Since
the function $\Delta(z)$ is holomorphic in the $G,$ arbitrary
point $\lambda$ form the $\mathcal{D}_0(T)$ is isolate in
$\mathcal{D}_0(T).$ Then the point $\lambda_0$ is isolate in the
spectrum $\sigma(T_1)\cup\sigma(T_2)\cup\mathcal{D}_0(T)$ of the
operator $T$ since $(\sigma(T_1)\cup\sigma(T_2))
\cap\mathcal{D}_0(T)=\emptyset.$ Thus, by definition of a discrete
spectrum, we obtain $\lambda_0\in\sigma_{disc}(T),$ i.e.
$\mathcal{D}_0(T)\subset\sigma_{disc}(T).$

Theorems \ref{EshkArz_thm4.6.} and \ref{EshkArz_thm5.5.} implies

\begin{thm}\label{EshkArz_thm5.6.}
$$\sigma_{ess}(T_1+T_2)=
\{0\}\cup\left(\bigcup\limits_{i=1}^n
Essran(h_i)\right)\cup\left(\bigcup\limits_{j=1}^m
Essran(p_j)\right).$$ \end{thm} \

\begin{exam}\label{EshArz_exam5.7.} Let $h_i(y)\equiv a_i \in\mathbb{R}\setminus
\{0\}, \,\ i\in\{1, ... , n\}$ and $p_j(x)\equiv
b_j\in\mathbb{R}\setminus \{0\}, \,\ j\in\{1, ... , m\}$ for the
kernels $k_1(x,s,y)$ and $k_2(x,t,y)$ of PIO $T_1$ and $T_2.$
\end{exam} Then by theorem \ref{EshkArz_thm2.3.} and
\ref{EshkArz_thm2.8.}, we obtain
$$\sigma(T_1)=\{0, a_1, ... , a_n\}, \,\,\ \sigma(T_2)=\{0, b_1, ... , b_m\}.$$
By theorem \ref{EshkArz_thm5.6.}, we have
$$\sigma_{ess}(T_1+T_2)=\{0, a_1, ...
, a_n, b_1, ... , b_m\}.$$ We obtain from \eqref{EshkArz_eq14}:
$$F_{k,j}(x,y;\lambda)=\frac{\lambda a_j}{(a_j-\lambda)(b_k-\lambda)}\varphi_j(x)\psi_k(y), \,\ \lambda\in\sigma_{ess}(T_1+T_2),$$ and
$$B_{k,j}(s,t)=b_k\overline{\varphi_j(s)}\cdot\overline{\psi_k(t)}.$$
Then the homogeneous Fredholm integral equation
\eqref{EshkArz_eq15} becomes the following form
\begin{equation}\label{EshkArz_eq17}
f(x,y)-\sum\limits_{j=1}^n\sum\limits_{k=1}^m\frac{a_jb_k}{(a_j-\lambda)(b_k-\lambda)}
\int_a^b\int_c^d\varphi_j(x)\psi_k(y)\overline{\varphi_j(s)}\cdot\overline{\psi_k(t)}d\mu_1(s)d\mu_2(t)=0.
\end{equation}

However, using the property of Fredholm integral equations with a
degenerate kernel, we obtain
$$\frac{a_jb_k}{(a_j-\lambda)(b_k-\lambda)}=1, \,\ j\in\{1,
... , n\}, \,\ k\in\{1, ... , m\}, \,\
\lambda\overline{\in}\sigma_{ess}(T_1+T_2).$$ It means that
$$a_jb_k=(a_j-\lambda)(b_k-\lambda), \,\ \lambda\overline{\in}\{0\}\cup\{a_i\}\cup\{b_l\},$$
i.e. the integral equation \eqref{EshkArz_eq17} has nontrivial
solution if and only if
$$\lambda=a_j+b_k\overline{\in} {\{0\}}\cup\{a_i\}\cup\{b_l\}.$$
Set
$$\Lambda=\{\lambda: \lambda=a_j+b_k\overline{\in}\{0, a_1, ...
, a_n, b_1, ... , b_m\}, \,\ j=\overline{1,n}, \,\
k=\overline{1,m}\}.$$

Therefore according to theorem \ref{EshkArz_thm5.5.}, we obtain
$\sigma_{disc}(T_1+T_2)=\Lambda$ and theorem \ref{EshkArz_thm4.6.}
implies
\begin{equation}\label{EshkArz_eq18}
\sigma(T_1+T_2)=\{\lambda: \lambda=a+b, \,\ a\in\sigma(T_1), \,\
b\in\sigma(T_2)\}. \end{equation}

It should be noted, the equality \eqref{EshkArz_eq18} was proved
for the PIO $T$ in $L_p, p\geq1$ with more general kernels
$k_1(x,s,y)=k_1(x,s),$ $k_2(x,t,y)=k_2(t,y)$ in the paper
\cite{kal4}, \cite{app2}.

\begin{rem}\label{EshArz_rem5.8.} It is known, that the discrete
spectrum $\sigma_{disc}(K)$ for each self-adjoint Fredholm
integral operator $K$ with a degenerate kernel is finite (since
$\sigma_{disc}(K)$ is the set of all eigenvalues of $K$ different
from zero). The following question is arisen: Does this property
hold for the PIO $T=T_1+T_2$ with a degenerate kernels (3) and (4)
This question is still an open problem.  \end{rem}

\end{document}